\documentclass[11pt]{article}
\usepackage{amsfonts, amsmath, amssymb, amsthm, float, fancyhdr}
\usepackage[pdftex]{graphicx}
\usepackage[a4paper,vmargin={20mm,20mm},hmargin={20mm,20mm}]{geometry}
\usepackage{enumerate}

\everymath{\displaystyle}
\DeclareGraphicsExtensions{.pdf,.png,.jpg}

\begin{document}
\begin{center}
{\huge{Permutations of counters on a table}}

\bigskip
Samuel Korsky

\bigskip
November 11, 2021

\bigskip
{\bf{Abstract}}
%\smallskip
\begin{changemargin}{2cm}{2cm} 
We consider a game in which a blindfolded player attempts to set $n$ counters lying on the vertices of a rotating regular $n$-gon table simultaneously to $0$. When the counters count$\pmod{m}$ we simplify the argument of Bar Yehuda, Etzion, and Moran (1993) showing that the player can win if and only if $n = 1$, $m = 1$, or $(n, m) = (p^a, p^b)$ for some prime $p$ and $a, b \in \mathbb{N}$.  We broadly generalize the result to the setting where the counters can be permuted by any element of a subset of the symmetric group $S \subseteq S_n$, with the original formulation corresponding to $S = \mathbb{Z}_n$ (rotations of the table).  
\end{changemargin} 
\end{center}
\section{Introduction}
We start with a well-known brainteaser which can be found in literature dating back to 1980 [1, 2, 3], and which was more recently popularized in 2019 on the famous prediction site fivethirtyeight.com [4]: 

\bigskip
\noindent
{\it{Four coins lie on the corners of a square table, some heads-up and some tails-up (they may all have the same orientation). Each turn, a blindfolded player can flip some of the coins, after which the table is rotated arbitrarily. If the player's goal is to at any time have all coins heads-up simultaneously, does he have a strategy that guarantees victory in a finite number of turns?}}

\bigskip
\noindent
For the simple case above, there is indeed a strategy that wins within $15$ turns. In particular, label the positions of the table $1, 2, 3, 4$, with these positions fixed from the perspective of the player. Then a {\it{move}} (performed once per turn) will consist of a vector in $\mathbb{Z}_2^4$, with $0$ denoting leaving the coin in the corresponding position as is and $1$ denoting flipping the coin in that position. E.g., the vector $(0, 0, 0, 1)$ denotes only flipping the coin in position $4$. The player's strategy should then be the following sequence of $15$ moves:

$$ (1, 1, 1, 1), (0, 1, 0, 1), (1, 1, 1, 1), (0, 0, 1, 1), (1, 1, 1, 1), (0, 1, 0, 1), (1, 1, 1, 1), (0, 0, 0, 1),$$
$$ (1, 1, 1, 1), (0, 1, 0, 1), (1, 1, 1, 1), (0, 0, 1, 1), (1, 1, 1, 1), (0, 1, 0, 1), (1, 1, 1, 1) $$

\bigskip
\noindent
It is easy to show by case work that these moves guarantee a win for the player.

\bigskip
\noindent
One naturally asks the question, what if instead of four coins there were $n$ coins? Furthermore, viewing a coin as a counter counting$\pmod{2}$ (to which the player adds either $0$ or $1$ each turn), what if instead the counters counted$\pmod{m}$, with the player adding one of $0, 1, \dots m-1$ to each counter each turn? This problem was considered by Bar Yehuda, Etzion, and Moran in [1], who showed:

\bigskip
\noindent
{\bf{Theorem 1.1}}. The player can win if and only if $n = 1$, $m = 1$, or $(n, m) = (p^a, p^b)$ for some prime $p$ and $a, b \in \mathbb{N}$.

\bigskip
\noindent
We independently derive this result, and simplify their argument by providing a clever explicit construction of a winning set of moves.  We also generalize the problem as follows: instead of the table simply rotating, one can imagine that the table permutes the counters in positions $1, 2, \dots, n$ based on elements from some subset $S \subseteq S_n$, where $S_n$ denotes the symmetric group and where $S$ contains the identity.\footnote{We may assume without loss of generality that $S$ contains the identity, since the player can pretend a certain permutation $t \in S$ happens every turn by default and then is followed by a permutation from the set $t^{-1} \cdot S$, which contains the identity.} Denote the parameters of this game by the ordered pair $(S, m)$, so that $(\mathbb{Z}_n, m)$ represents the setting in Theorem 1.1. Additionally, let $G \le S_n$ be the subgroup of $S_n$ generated by $S$. Our main result is:

\bigskip
\noindent
{\bf{Theorem 1.2}}. The player can win the $(S, m)$-game if and only if $|G| = 1$, $m = 1$, or $(|G|, m) = (p^a, p^b)$ for some prime $p$ and $a, b \in \mathbb{N}$.

\bigskip
\noindent
Our paper is divided into five parts, where the first four parts consist of a simplification of the proof of Theorem 1.1 in [1].  In Section 2, we show that if $(n, m) = (p, q)$ for distinct primes $p, q$ then the player cannot win. In Section 3, we show if the player cannot win when $(n, m) = (a, b)$ then the player also cannot win for $(n, m) = (a, bk)$ or $(n, m) = (ak, b)$ for any $k \in \mathbb{N}$. In Section 4, we constructively show that the player can win if $(n, m) = (p^a, p)$ for a prime $p$ and any $a \in \mathbb{N}$. In Section 5, we extend the construction to the case where $(n, m) = (p^a, p^b)$ where $b > 1$. Finally, drawing on the methods used in Sections 2 through 5, in Section 6 we prove Theorem 1.2 in full generality. 

\section{The $(n, m) = (p, q)$ case}
We start by extending the notation in the introduction. Instead of denoting ``moves" by vectors in $\mathbb{Z}_2^4$, we will now use vectors in $\mathbb{Z}_m^n$. E.g., the vector $(2, 0, 3)$ will denote adding $2$ to the counter in position $1$ and adding $3$ to the counter in position $3$. We will also use vectors in $\mathbb{Z}_m^n$ to describe the configuration of counters. Furthermore, call a configuration of counters {\it{homogenous}} if each counter on the table shows the same number, and non-homogenous otherwise. 

\bigskip
\noindent
{\bf{Lemma 2.1}}. The player cannot win if $(n, m) = (p, q)$ for distinct primes $p, q$. 

\bigskip
\noindent
{\it{Proof.}} When $(n, m) = (p, q)$ for distinct primes $p, q$, we will show that for any non-homogenous configuration, there is no move guaranteed to make the configuration homogenous following an arbitrary rotation of the table. 

\bigskip
\noindent
Indeed, suppose the configuration on the table prior to a rotation was $(x_1, x_2, \dots, x_p)$ and consider any move $(y_1, y_2, \dots, y_p)$. For this move to guarantee that the configuration of the table afterwards was homogenous, the following equalities would have to hold simultaneously:
\begin{align*}
x_1 + y_1 = x_2 + y_2 = &\dots = x_p + y_p\pmod{q}\\
x_p + y_1 = x_1 + y_2 = &\dots = x_{p-1} + y_p\pmod{q}\\
&\vdots\\
x_2 + y_1 = x_3 + y_2 = &\dots = x_1 + y_p\pmod{q}.
\end{align*}
This implies that
$$ x_1 - x_2 =x_2 - x_3 = \dots = x_p - x_1 = y_2 - y_1\pmod{q}. $$
Thus
$$ p(x_1 - x_2) = (x_1 - x_2) + (x_2 - x_3) + \dots + (x_p - x_1) = 0\pmod{q}$$
and so $x_1 = x_2\pmod{q}$. Similarly we obtain $x_1 = x_2 = \dots = x_p\pmod{q}$ so indeed such a move is only possible if the configuration of the table was already homogenous.

\bigskip
\noindent
Therefore if the starting configuration is non-homogenous, the player can never force the configuration to be homogenous and so cannot win. $\blacksquare$

\section{The $(n, m) = (ak, b)$ and $(n, m) = (a, bk)$ cases}
We handle each case separately:

\bigskip
\noindent
{\bf{Lemma 3.1}}. If the player cannot win when $(n, m) = (a, b)$, then the player cannot win when $(n, m) = (ak, b)$ for any $k \in \mathbb{N}$. 

\bigskip
\noindent
{\it{Proof.}} Suppose that the player cannot win if $(n, m) = (a, b)$ for some $(a, b) \in \mathbb{N}^2$, and consider the case where $(n, m) = (ak, b)$ for some $k \in \mathbb{N}$. Suppose for the sake of contradiction that the player had a sequence of moves $y_1, y_2, \dots, y_N$ for some $N \in \mathbb{N}$ that guaranteed a win, where $y_i = (y_{i,1}, y_{i,2}, \dots, y_{i,ak})$ for all $i$. Now let $y'_i =(y_{i,k}, y_{i,2k}, \dots, y_{i,ak})$ for all $i$. We must have that the sequence of moves $y'_1, y'_2, \dots, y'_N$ wins for $(n, m) = (a, b)$, contradiction. $\blacksquare$

\bigskip
\noindent
{\bf{Lemma 3.2}}. If the player cannot win when $(n, m) = (a, b)$, then the player cannot win when $(n, m) = (a, bk)$ for any $k \in \mathbb{N}$.

\bigskip
\noindent
{\it{Proof.}} Suppose that the player cannot win if $(n, m) = (a, b)$ for some $(a, b) \in \mathbb{N}^2$, and consider the case where $(n, m) = (a, bk)$ for some $k \in \mathbb{N}$. Suppose for the sake of contradiction that the player had a sequence of moves $y_1, y_2, \dots, y_N$ for some $N \in \mathbb{N}$ that guaranteed a win, where $y_i =(y_{i,1}, y_{i,2}, \dots, y_{i,a})$ for all $i$. Note that $y_{i,j} \in \mathbb{Z}_{bk}$ for all $i, j$, and define the homomorphism $\phi: \mathbb{Z}_{bk} \rightarrow \mathbb{Z}_b$ where $\phi(x) = x\pmod{b}$. Now let $y'_i = (\phi(y_{i,1}), \phi(y_{i,2}), \dots, \phi(y_{i,a}))$ for all $i$. We must have that the sequence of moves $y'_1, y'_2, \dots, y'_N$ wins for $(n, m) = (a, b)$, contradiction. $\blacksquare$

\bigskip
\noindent
{\bf{Corollary 3.3}}. The combination of Lemmas 2.1, 3.1, and 3.2 immediately imply the ``only if" direction in Theorem 1.1. 

\section{The $(n, m) = (p^a, p)$ case}
Here we prove the following lemma constructively:

\bigskip
\noindent
{\bf{Lemma 4.1}}. The player can win if $(n, m) = (p^a, p)$ for some prime $p$ and $a \in \mathbb{N}$. 

\bigskip
\noindent
{\it{Proof.}} Let $x_{i,j} = \binom{i}{j}\pmod{p}$ where $x_{i,j} \in \mathbb{Z}_p$ for all $i,j$. Furthermore, let $x_j = (x_{0,j}, x_{1,j}, \dots, x_{p^a-1,j})$ for all $j \in \{0, 1, \dots, p^a - 1\}$. For all $i \in \{1, 2, \dots, p^{p^a} - 1\}$, let $y_i = x_{v_p(i)}$ where $v_p$ denotes $p$-adic valuation. We claim that the sequence of moves $y_1, y_2, \dots, y_{p^{p^a} - 1}$ wins.

\bigskip
\noindent
Since the matrix $[x_0^T\ x_1^T\ \dots\ x_{p^a-1}^T]$ is lower triangular and its main diagonal is identically $1$, its determinant is $1 \ne 0$ and so the moves $x_0, x_1, \dots, x_{p-1}$ form a basis over $\mathbb{Z}_p^{p^a}$. Therefore, we can write the starting configuration $s$ as $s = c_0x_0 + c_1x_1 + \dots c_{p^a-1}x_{p^a-1}$ for some $c_0, c_1, \dots, c_{p^a-1} \in \mathbb{Z}_p$. 

\bigskip
\noindent
The following intermediate claim is then the key to the proof:

\bigskip
\noindent
{\bf{Claim 4.2}}. For any $j$, let $x'_j$ be a cyclic permutation of $x_j$. Then $x_j - x'_j = e_0x_0 + e_1x_1 + \dots + e_{j-1}x_{j-1}$ for some $e_0, e_1, \dots, e_{j-1} \in \mathbb{Z}_p$. 

\bigskip
\noindent
{\it{Proof.}} We proceed by induction on $j$. When $j = 0$ the result is trivial, since $x_j = x'_j$. Now suppose $j > 0$ and let $x^{(k)}_j = (x_{k,j}, x_{k+1,j}, \dots, x_{k-1,j})$ for all $k \in \{0, 1, \dots, p^a - 1\}$, so that $x^{(0)}_j = x_j$. Repeatedly utilizing the fact that $\binom{i}{j} - \binom{i-1}{j} = \binom{i-1}{j-1}$ and $\binom{p^a}{j} = \binom{0}{j} = 0\pmod{p}$ we have that, working in $\mathbb{Z}_p^{p^a}$,
\begin{align*} 
x^{(k+1)}_j - x^{(k)}_j &= (x_{k+1,j}, x_{k+2,j}, \dots, x_{k,j}) - (x_{k,j}, x_{k+1,j}, \dots, x_{k-1,j})\\
&= (x_{k,j-1}, x_{k+1,j-1}, \dots, x_{k-1,j-1})\\
&= x^{(k)}_{j-1}
\end{align*}
Therefore letting $x'_j = x^{(k)}_j$ for some $k$ we have
\begin{align*}
x'_j - x_j &= \left(x^{(k)}_j - x^{(k-1)}_j\right) +\left (x^{(k-1)}_j - x^{(k-2)}_j\right) + \dots + \left(x^{(1)}_j - x^{0)}_j\right)\\
&= x^{(k-1)}_{j-1} + x^{(k-2)}_{j-1} + \dots + x^{(0)}_{j-1}
\end{align*}
But by the inductive hypothesis we know that each $x^{(i)}_{j-1}$ can be written as a linear combination of $x_0, x_1, \dots, x_{j-1}$ (with a coefficient of $1$ behind $x_{j-1}$), which completes the proof. $\blacksquare$

\bigskip
\noindent
Returning to the proof of Lemma 4.1, recall that we wrote $s = c_0x_0 + c_1x_1 + \dots c_{p^a-1}x_{p^a-1}$ for some $c_0, c_1, \dots, c_{p^a-1} \in \mathbb{Z}_p$. For any such starting configuration $s$, let $f(s)$ denote that largest $i$ such that $c_i \ne 0$. We will show by induction on $f(s)$ that the sequence of moves $y_1, y_2, \dots, y_{p^{f(s)+1} - 1}$ wins. Note that the base case where $c_i = 0$ for all $i$ is trivial, since the player immediately wins. 

\bigskip
\noindent
The main idea is that as a consequence of Claim 4.2, every time the table rotates and we rewrite the new configuration as a linear combination of $x_0, x_1, \dots, x_{p^a-1}$, the coefficient behind $x_{f(s)}$ is invariant. Specifically, suppose $c_{f(s)} = c \ne 0$. Let $s'$ be the configuration on the table after $(p - c)p^{f(s)}$ moves, and write $s' = c'_0x_0 + c'_1x_1 + \dots c'_{p^a-1}x_{p^a-1}$ for some $c'_0, c'_1, \dots, c'_{p^a-1} \in \mathbb{Z}_p$. By Claim 4.2, each move $y_i$ with $v_p(i) < f(s)$ would not affect any of the coefficients behind $x_j$ for any $j \ge f(s)$ and each of the $p - c$ moves with $v_p(i) = f(s)$ would increase the coefficient behind $x_{f(s)}$ by $1$ regardless of how the table rotates between moves, so that $c'_{f(s)} = c_{f(s)} + p - c = 0\pmod{p}$. Therefore after $(p - c)p^{f(s)}$ moves we are in a configuration with $f(s') < f(s)$ and since the next $p^{f(s') + 1} - 1$ moves are copies of the first $p^{f(s') + 1} - 1$ moves, we are done by induction. $\blacksquare$

\section{The $(n, m) = (p^a, p^b)$ case}
Finally, we expand upon our construction in Section 4 to prove Theorem 1.1 in its entirety. Specifically, we will show by induction on $b$ that there is a sequence of $p^{bp^a} - 1$ moves that wins. The base case of $b = 1$ follows from the proof of Lemma 4.1. 

\bigskip
\noindent
Now, suppose $x_1, x_2, \dots, x_{p^{(b - 1)p^a} - 1}$ is a sequence of moves that wins in the $(n, m) = (p^a, p^{b-1})$ case. Additionally, let $y_1, y_2, \dots, y_{p^{p^a} - 1}$ be the sequence of moves that wins in the $(n, m) = (p^a, p)$ case as in the proof of Lemma 4.1. Note that $x_i \in \mathbb{Z}_{p^{b-1}}^{p^a}$ for all $i$ and $y_j \in \mathbb{Z}_p^{p^a}$ for all $j$, but interpret each of these vectors as vectors in $\mathbb{Z}_{p^b}^{p^a}$ (through the identity homomorphism). Define the homomorphism $\varphi: \mathbb{Z}_{p^b}^{p^a} \rightarrow \mathbb{Z}_{p^{b-1}}^{p^a}$ where $\varphi(x) = x\pmod{p^{b-1}}$ element-wise, and let
$$ z_i = \begin{cases} px_{\varphi(i)} &\mbox{if } i \ne 0 \\ y_{ip^{(1-b)p^a}} & \mbox{if } i = 0 \end{cases}\pmod{p^{(b -1)p^a}}$$
for all $i \in \{1, 2, \dots, p^{bp^a} - 1\}$, where $px_{\varphi(i)}$ denotes element-wise multiplication by $p$. We claim that $z_1, z_2, \dots, z_{p^{bp^a} - 1}$ is a sequence of moves that wins in the $(n, m) = (p^a, p^b)$ case. 

\bigskip
\noindent
If we consider the configuration of the table $\pmod{p}$, then up to rotation none of the moves of the form $px_{\varphi(i)}$ affect the configuration. Therefore we know by the definition of the sequence $y_1, y_2, \dots, y_{p^{p^a} - 1}$ that after some move of the form $y_j$, every counter will be $0\pmod{p}$. It is then clear that after this move $y_j$, the following sequence of $p^{(b -1)p^a} - 1$ moves $px_1, px_2, \dots, px_{p^{(b-1)p^a} - 1}$ will win the game. This completes the proof of Theorem 1.1. $\blacksquare$

\bigskip
\noindent
We will also prove that our constructions are optimal in terms of number of moves necessary to win. 

\bigskip
\noindent
{\bf{Theorem 5.1}}. If the player can win for a given $(n, m)$, then he can guarantee a win in no less than $m^n - 1$ moves. 

\bigskip
\noindent
{\it{Proof}}. Consider any sequence of moves $y_1, y_2, \dots, y_N \in \mathbb{Z}_m^n$ with $N < m^n - 1$. Let $z_k = \sum_{i=1}^{k}y_i$ for all $k \in \{1, 2, \dots, N\}$ and suppose that the table did not rotate at all after any of the moves. Then this sequence of moves would only win if the starting configuration was equal to $0$ or $-z_k$ (with each element reduced$\pmod{m}$) for some $k$. But there are $N + 1 < m^n$ such winning configurations and $m^n$ possible starting configurations, so there are a non-zero number of starting configurations for which this sequence of moves never wins. This implies the desired result. $\blacksquare$

\section{Generalizing to groups}
Consider any subset of permutations $S \subseteq S_n$ with the identity ${\bf{1}} \in S$ and suppose that instead of simply rotating each turn, the counters on the table can be acted on by any permutation in $S$.  As in the introduction, let $G \le S_n$ be the subgroup of $S_n$ generated by $S$.  Here we will prove Theorem 1.2, which states that the player can win the $(S, m)$-game if and only if $|G| = 1$, $m = 1$, or $(|G|, m) = (p^a, p^b)$ for some prime $p$ and $a, b \in \mathbb{N}$.

\bigskip
\noindent
We break the proof into two parts: first we shall show the ``only if" direction, and then we shall show the ``if" direction when $b = 1$.  The case where $b > 1$ will then immediately follow from the logic in the proof of Theorem 1.1 in Section 5. 

\bigskip
\noindent
{\bf{Lemma 6.1}}. The player can win the $(S, m)$-game only if $|G| = 1$, $m = 1$, or $(|G|, m) = (p^a, p^b)$ for some prime $p$ and $a, b \in \mathbb{N}$.

\bigskip
\noindent
{\it{Proof}}. We mimic the proof of Lemma 2.1.  Suppose that there exist distinct primes $p$ and $q$ with $v_p(|G|), v_q(m) > 0$ and without loss of generality assume $m = q$. By Cauchy's Theorem there must exist some $c \in G$ with order $p$.  Let $g(x)$ denote the position of the counter currently at position $x$ after the permutation $g \in G$ is applied to the counters on the table.  Call a configuration of the table {\it{semi-homogenous}} if the counters in positions $g(1), gc(1),gc^2(1), \dots, gc^{p-1}(1)$ show the same number for all $g \in G$. We will show that for any non-semi-homogenous configuration on the table, there is no move guaranteed to make the configuration semi-homogenous following an arbitrary permutation of the table.

\bigskip
\noindent
Indeed, suppose the configuration on the table prior to a permutation was $(x_1, x_2, \dots, x_n)$ and consider any move $(y_1, y_2, \dots, y_n)$.  For convenience, let $x_g$ denote $x_{g(1)}$ for all $g \in G$ and define $y_g$ similarly, and let us work in $\mathbb{Z}_q$.  Additionally, let $T = \{s^{-1} | s \in S\}$ be the set of inverses of elements in $S$.  For this move to guarantee that the configuration of the table was semi-homogenous after the move, the following $|S||G|$ strings of equalities would have to hold simultaneously:
$$ x_{tg} + y_g = x_{tgc} + y_{gc} = \dots = x_{tgc^{p-1}} + y_{gc^{p-1}} $$
for all $t \in T$ and $g \in G$.  Now, fix a specific $d \in G$ and write $dcd^{-1} = \prod_{i = 1}^{k}t_i$ for some $k \in \mathbb{N}$ and $t_1, t_2, \dots, t_k \in T$ (this representation is guaranteed to exist since $T$ generates $G$).  Using the first equality in the string with $g = d$ and $t \in \{{\bf{1}}, t_k\}$, we obtain $x_d - x_{dc} = x_{t_kd} - x_{t_kdc} =y_{dc} - y_d$.  Using the first equality again with $g = t_kd$ and $t \in \{{\bf{1}}, t_{k-1}\}$ we obtain $x_{t_kd}  - x_{t_kdc} = x_{t_{k-1}t_kd} - x_{t_{k-1}t_kdc} =y_{t_kdc} - y_{t_kd}$.  Combining equalities we obtain $x _d - x_{dc} = x_{t_{k-1}t_kd} - x_{t_{k-1}t_kdc}$. Continuing in this fashion we obtain 
$$ x_d - x_{dc} = x_{t_1t_2{\dots}t_kd} - x_{t_1t_2{\dots}t_kdc} = x_{dc} - x_{dc^2} $$
and repeating the argument we have
$$  x_d - x_{dc} = x_{dc} - x_{dc^2} = \dots =  x_{dc^{p-1}} - x_{d} $$
which holds for any $d \in G$.  Notice that these $p$ expressions sum to $0$, so since $p$ and $q$ are distinct each expression must equal $0$ and so the configuration $(x_1, x_2, \dots, x_n)$ must have been semi-homogenous to begin with. 

\bigskip
\noindent
Therefore if the starting configuration is not semi-homogenous, the player can never force the configuration to be semi-homogenous and so cannot win. $\blacksquare$

\bigskip
\noindent
Now we proceed to the second part of the proof of Theorem 1.2:

\bigskip
\noindent
{\bf{Lemma 6.2}}. The player can win the $(G, m)$-game if $(|G|, m) = (p^a, p)$ for some prime $p$ and $a \in \mathbb{N}$.

\bigskip
\noindent
{\it{Proof}}. We mimic the proof of Lemma 4.1.  Suppose there exist vectors $x_0, x_1, \dots, x_{n-1}$ that form a basis for $\mathbb{Z}_p^n$ and that have the property that $x_j - g \cdot x_j$ can be written as a linear combination of $x_1, x_2, \dots, x_{j-1}$ for all $j \in \{0, 1, \dots, n-1\}$ and all $g \in G$, where $g \cdot x$ represents the permutation of the coordinates of $x$ corresponding to $g \in S_n$.  For all $i \in \{1, 2, \dots, p^{n} - 1\}$, let $y_i = x_{v_p(i)}$ where $v_p$ denotes $p$-adic valuation.  Then by the same reasoning as from the proof of Lemma 4.1, the sequence of moves $y_1, y_2, \dots, y_{p^{n} - 1}$ wins.

\bigskip
\noindent
Thus it suffices to show that such a basis $x_0, x_1, \dots, x_{n-1}$ exists. Consider the orbits of $G$ on $\mathbb{Z}_p^n$, and suppose we partition $\mathbb{Z}_p^n$ into these orbits.  Since $|G| = p^a$, the Orbit-Stabilizer Theorem implies that each orbit has size $p^k$ for some $k \in \{0, 1, 2, \dots, a\}$.  Since the number of vectors in $\mathbb{Z}_p^n$ is $p^n$, the number of orbits of size $1$ must be divisible by $p$.  But note that $0$ has an orbit of size $1$, so there exists some nonzero $x_0 \in  \mathbb{Z}_p^n$ fixed by $G$.  

\bigskip
\noindent
We can repeat the argument on the quotient space $\mathbb{Z}_p^n/\langle x_0 \rangle$ to find some nonzero $x_1 \in \mathbb{Z}_p^n/\langle x_0 \rangle$ that is fixed by $G$. Continuing in this fashion we obtain a nested sequence of subspaces 
$$ \langle x_0 \rangle < \langle x_0, x_1 \rangle < \dots < \langle x_0, x_1, \dots, x_{n-1} \rangle $$
each of which is fixed by $G$, and it is clear that the basis $x_0, x_1, \dots, x_{n-1}$ satisfies the desired condition. This completes the proof. $\blacksquare$

\bigskip
\noindent
Given a sequence of moves that wins the $(G, p)$-game, we can then induct on $b$ as in Section 5 to construct a sequence of moves that wins the $(G, p^b)$-game for any $b \in \mathbb{N}$. Furthermore, if the player can win the $(G, m)$-game then since $S \subseteq G$ he can use the same sequence of moves to win the $(S, m)$-game, so the combination of Lemmas 6.1 and 6.2 imply Theorem 1.2, as desired. $\blacksquare$

\section{Acknowledgements}
The author would like to thank Dhroova Aiylam and Alexander Katz for their helpful discussions.

\end{document}